%
%
%
%
\def\references{sperner.ref}        
\baselineskip = 15.5pt plus .5pt minus .5pt  
%
\font\bigbf=cmbx10   scaled \magstep1     
\font\smallcaps = cmcsc10                 
%
\ifx\fonts\cmfonts
\font\ninerm=cmr9
\font\ninei=cmmi9
\font\ninesy=cmsy9
\font\ninebf=cmbx9
\font\ninett=cmtt9
\font\nineit=cmti9
\font\ninesl=cmsl9
\else
\font\ninerm=amr9
\font\ninei=ammi9
\font\ninesy=amsy9
\font\ninebf=ambx9
\font\ninett=amtt9
\font\nineit=amti9
\font\ninesl=amsl9
\fi
\skewchar\ninei='177
\skewchar\ninesy='60
\skewchar\ninett=-1
\newskip\tglue
\def\ninepoint{\def\rm{\fam0\ninerm}
       \textfont0=\ninerm
       \textfont1=\ninei
       \textfont2=\ninesy
       \textfont\itfam=\nineit \def\it{\fam\itfam\nineit}%
       \textfont\slfam=\ninesl \def\sl{\fam\slfam\ninesl}%
       \textfont\ttfam=\ninett \def\tt{\fam\ttfam\ninett}%
       \textfont\bffam=\ninebf \def\bf{\fam\bffam\ninebf}%
       \tt\tglue=.5em plus.25em minus .15em
       \normalbaselineskip=11pt
       \setbox\strutbox=\hbox{\vrule height8pt depth3pt width0pt}%
       \let\sc=\sevenrm \let\big=\ninebig \normalbaselines\rm}%
\input amssym
\input \references  
\def\fonts{cmfonts}
%
\def\section#1{\vskip0pt plus .1\vsize
    \penalty-250\vskip0pt plus-.1\vsize\bigskip
    \noindent{\bf #1}\nobreak\message{#1}}

\def\abstract#1{\bigskip\centerline{\hbox{
       \vbox{\hsize=4.75truein%
           {\ninepoint \noindent ABSTRACT.\enspace#1}
            }}}}
\def\keywords#1{\bigskip\centerline{\hbox{
       \vbox{\hsize=4.85truein%
           {\ninepoint \noindent KEYWORDS.\enspace#1}
            }}}}
\def\AmsClass#1#2{\bigskip\centerline{\hbox{
       \vbox{\hsize=4.75truein%
          {\ninepoint \noindent
       AMS #1 Subject Classification Numbers:\enspace#2}
            }}}}
\def\th #1 #2: #3\par\par{\medbreak{\bf#1 #2:
\enspace}{\sl#3\par}\par\medbreak}
\def\co #1 #2: #3\par\par{\medbreak{\bf#1 #2:
\enspace}{\sl#3\par}\par\medbreak}
\def\le #1 #2: #3\par\par{\medbreak{\bf #1 #2:
\enspace}{\sl #3\par}\par\medbreak}
\def\rem #1 #2. #3\par{\medbreak{\bf #1 #2.
\enspace}{#3}\par\medbreak}
\def\proof{{\bf Proof}.\enspace}
\def\diam{\hbox{\rm diam}}

\def\sqr#1#2{{\vcenter{\hrule height.#2pt
      \hbox{\vrule width.#2pt height#1pt \kern#1pt
         \vrule width.#2pt}
       \hrule height.#2pt}}}

\def\boxit#1{\vbox{\hrule \hbox{\vrule \kern2pt
                 \vbox{\kern2pt#1\kern2pt}\kern2pt\vrule}\hrule}}
%
%
\newdimen\refindent\newdimen\plusindent
\newdimen\refskip\newdimen\tempindent
\newdimen\extraindent
\newcount\refcount
\newwrite\reffile
\def\beginref{}
\def\endref{}
\def\ifundefined#1{\expandafter\ifx\csname#1\endcsname\relax}
\def\referto[#1]{\ifundefined{#1}[?]\else[\csname#1\endcsname]\fi}
%
%
\refcount=0
\def\ref#1:#2.-#3[#4]{\ninepoint 
\advance\refcount by 1
\setbox0=\hbox{[\number\refcount]}\refindent=\wd0
\plusindent=\refskip\extraindent=\refskip
\advance\plusindent by -\refindent\tempindent=\parindent %
\parindent=0pt\par\hangindent\extraindent %
 [\number\refcount]\hskip\plusindent #1:{\sl#2},#3
\parindent=\tempindent
}
\refskip=\parindent

\def\int{\hbox{\rm int}}

\def\c#1{{\cal #1}}

\def\beth{\hbox{\tt [}}

%
\def\specialheadlines{
        \headline={\vbox{\line{
                     \ifnum\pageno<2 \ffolio
                     \else\rightheadline
                     \fi}\bigskip
                     }}}
        \def\ffolio{\hfil}
        \def\rightheadline{\hfil
                    {\smallcaps An Infinitary Version of
                    Sperner's Lemma}\hfil
                    \hbox{\rm\folio}}
\nopagenumbers
\specialheadlines
%
%
\null
\vskip2truecm
  \centerline{\bigbf An Infinitary Version of Sperner's Lemma}
\vskip.75truecm
 \centerline{\sl Aarno Hohti, University of Helsinki, Finland}
%
\vskip.5truecm
\abstract{We prove an extension of the well-known
          combinatorial-topological lemma of E.\ Sperner
          \referto[Spernera] to the
          case of infinite-dimensional cubes. It is obtained as a corollary
          to an infinitary extension of the Lebesgue Covering
          Dimension Theorem.}


\smallskip
\AmsClass{1991}{Primary 57N20, Secondary 55M20, 54F45.}

\vskip.5truecm
%
\section{1. Introduction.}
The well-known lemma of E.\ Sperner on colourings of vertices of
$n$-simplices is one of the cornerstones of simplicial algebraic
topology, directly related to the invariance of homology
groups under simplicial subdivisions. It was used by Sperner
\referto[Spernera]
to give a short proof that the topological dimension of an $n$-cell
equals $n$. In fact, it can be considered a combinatorial (discrete) version of
Brouwer's Fixed Point Theorem, obtained from Sperner's Lemma by a
simple argument using the compactness of the cubes $I^n$, as shown
in [8]. There
are several versions and
extensions of this lemma
(\referto[Browna],
\referto[Fana],
\referto[Kuhna],
\referto[Mania],
\referto[Spernerb],
\referto[Tuckera]) and there are
versions for matroids
(\referto[Krynskia],
\referto[Lindstroma],
\referto[Lovasza]). The well-known extension of
Brouwer's Fixed Point Theorem to infinite-dimensional (Banach) spaces,
Schauder's Fixed Point Theorem, is stated for {\it compact} convex
subsets and is based on the finite-dimensional theorem. Indeed,
there is no ``truly'' infinitary Sperner lemma in the literature, and
this is connected with the facts that Brouwer's Fixed Point Theorem
fails in general (see e.\ g.\ \referto[Benyaminia]) and the ordinary homology
groups $H_n(U,\Bbb Z)$ vanish for unit spheres $U$ of infinite-dimensional
Banach spaces. However, we give here a natural extension of Sperner's
Lemma to colourings of cubical triangulations of
infinite-dimensional cubes. The problem
can be reduced to a combinatorial problem about colourings
$\varphi : [k]^{\omega}\to \{0,1\}^{\omega}$ (where $k$ is a
positive integer and $[k]$ denotes
the set $\{0,\ldots,k\}$) satisfying the {\it Sperner condition} that
$\varphi(\sigma)\neq \varphi(\sigma')$ whenever the distance of
$\sigma$ and $\sigma'$ in $[k]^{\omega}$ is maximal. For such
colourings $\varphi$ we prove that there is $\sigma\in [k]^{\omega}$
such that $\varphi(K_{\sigma})$ is infinite, where $K_{\sigma}$ is
the ``cube'' corresponding to $\sigma$. We also indicate why
this result is the best possible. The results are stated in this
paper for the unit cube $U_{\infty}$ of the Banach space $\ell_{\infty}$;
we consider $U_{\infty}$ as the $\omega$-dimensional (combinatorial!)
analogue of the finite-dimensional cubes $I^n$.

The Lebesgue Covering Dimension Theorem
(cf., e.\ g., \referto[Engelkinga]) says that open covers
of $I^n$ by sufficiently small sets have order at least $n + 1$.
This theorem cannot be extended to open covers of the cube $U_{\infty}$,
since by paracompactness every open cover has a point-finite refinement.
However, this theorem can be extended if only {\it uniform}
covers are considered, replacing the lack of compactness of $U_{\infty}$
by the condition of uniformity. In this paper,
this extension is proved first and then used to obtain Sperner's
Lemma as a direct corollary.
The sets $[k]^{\omega}$
correspond to regular or uniform cubical subdivisions of $U_{\infty}$,
and, by the same token,  Sperner's Lemma is not valid for
non-uniform subdivisions. The general question
(Stone \referto[Stonea], Isbell \referto[Isbella]) whether every uniform
cover of a Banach space has a point-finite uniform refinement was
answered in the negative independently by
Pelant \referto[Pelanta] and \v S\v cepin
\referto[Scepina]. It was handled in a different way
later by R\"odl \referto[Rodla], who by using results of
\referto[Erdosa] produced
a counter-example of minimal cardinality.

\section{2. Preliminaries.} In this section we develop the necessary
notation and background for our treatment of the infinitary
covering dimension and
Sperner's Lemma. Let us first state the classical version of this
result. In this paper the symbol $[n]$ denotes
the subset $\{0,\ldots,n\}$ of integers.
Let $\Delta$ be an $n$-simplex, and let
$\c{K} = \{\Delta_1,\ldots,\Delta_m\}$ be a subdivision of $\Delta$.
Let $\varphi : \c{K}^{(0)}\to [n]$ be a mapping (``colouring'') of the
vertices of the simplices in $\c{K}$.  If $\varphi$ satisfies
the condition
\medskip
\noindent
{\bf (Sperner Condition)}:
$\varphi(\Delta^{(0)}) = [n]$ (bijectivity) and
if $v\in \c{K}^{(0)}$ lies in an
$(n-1)$-face $F$ of $\Delta$, then $\varphi(v)\in\varphi(F^{(0)}),$
\medskip
\noindent
then $\varphi$ is called a Sperner colouring. (In some articles
the notion ``proper labelling'' is used.)
The classical Sperner Lemma asserts that
there is a simplex $\Delta_i\in \c{K}$ such that
$\varphi(\Delta_i^{(0)}) = [n]$:

\th Theorem 2.1(Sperner's Lemma \referto[Spernera]):
If $\c{K}$ is a subdivision of an
$n$-simplex $\Delta$ and $\varphi:\c{K}^{(0)}\to [n]$ is
a Sperner colouring, then there is a simplex $\Delta_i$ of
$\c{K}$ such that the vertices of $\Delta_i$ are coloured
by $\varphi$ with $n+1$ colours.

In the sequel we loosely call colourings satisfying a suitable version of
Sperner
Condition {\it Sperner colourings}. We notice here that
in this finite-dimensional case, Sperner Condition
can be replaced by the condition that the cover
$\c{C} = \{\varphi^{-1}(k): k\in [n]\}$ determined by the
colours of the vertices $v\in \Delta^{(0)}$ is a ``bounded'' cover
of $\c{K}^{(0)}$ in the following sense. The complex $\c{K}$
becomes a {\it uniform complex} in the sense of Isbell \referto[Isbella] if
for $p,q\in |\c{K}|$ (here $|\c{K}|$ denotes the underlying
set of $\c{K}$) we let the distance $d(p,q)$ be the
maximum difference of their barycentric coordinates. Then
$\c{C}$ is called bounded if $\hbox{\rm diam}_d(\varphi^{-1}(k)) < 1$
for all $k\in [n]$. Notice that boundedness implies
$\varphi\upharpoonright \Delta^{(0)}$ is injective: if
$v,w\in \Delta^{(0)}$ and $v\neq w$, then $d(v,w) = 1$ and hence
$\varphi(v)\neq\varphi(w)$. Thus, $\varphi\upharpoonright\Delta^{(0)}$ is
bijective.

It is natural to consider cubes $I^n = [0,1]^n$ instead of the
simplices $\Delta$. For a ``cubical'' version of Sperner's Lemma, see
\referto[Kuhna]. In the infinite-dimensional situation, the cubical form
becomes a natural one. Let $\c{K}$ be the cubical complex
consisting of the single cube $I^n$. (We consider cubical complexes given by
a set of {\it maximal} cubes such that any two intersecting cubes
meet in a common cubical face.)
Then $\c{K}^{(0)}$ corresponds to the set $\{0,1\}^n$. A regular (or
uniform) subdivision of $I^n$ of sidelength $1/m$ ($m\geq 1$)
is the set consisting of all
products (called cubes of the subdivision)
$$K_\sigma \;=\; [i_1/m,(i_1+1)/m]\times\cdots\times [i_n/m,(i_n+1)/m],$$
where $i_1,\ldots,i_n\in [m-1]$ and where
the $n$-tuple $\sigma = (i_1,\ldots,i_n)$ is called the {\it index}
of the associated cube $K_{\sigma}$. For a
regular subdivision $\c{L}$ of $I^n$, the vertex set
$\c{L}^{(0)}$ corresponds to $[m]^n$ for some $m$; the correspondence
is simply given by the map $f: [m]^n\to \c{L}^{(0)}$ for which
$$f(i_1,\ldots,i_n) = (i_1/m,\ldots,i_n/m).$$
Thus, the colourings of $\c{L}^{(0)}$ satisfying Sperner Condition lead
to colourings $\varphi :[m]^n\to \{0,1\}^n$ with the following property:
if $\sigma,\sigma'\in [m]^n$ and $\sigma(i) = 0$, $\sigma'(i) = m$
for some $i\in [n]$, then $\varphi(\sigma)\neq\varphi(\sigma')$.
The natural way to extend these colourings to the infinitary
situation is to consider  colourings
$\varphi: [m]^{\omega}\to \{0,1\}^{\omega}$, where $m\in \omega$. As above,
the sets $[m]^{\omega}$ correspond to regular cubical subdivisions
of the cube $[0,1]^{\omega}$, denoted here by $U_{\infty}$. The
topology of $U_{\infty}$ is given by the $\ell_{\infty}$-norm
defined by $||\sigma - \tau||_{\ell_{\infty}} =
\sup\{|\sigma(i) - \tau(i)|: i\in\omega\}.$ We denote the distance
of two elements $\sigma,\tau\in [m]^{\omega}$ simply by
$||\sigma - \tau||$.

\section{3. Regular subdivisions.} The classical Sperner Lemma
was extended to colourings of
cubical triangulations by Kuhn in
\referto[Kuhna]. It is easy to
show by using either the classical result that $\dim(I^n) = n$
or the results of Kuhn that for any colouring
$\varphi : [m]^n\to \{0,1\}^n$ satisfying Sperner Condition
there is $\sigma\in [m]^n$ such that $\varphi(K_{\sigma})$
contains at least $n+1$ colours. (The result proved by Kuhn
is stronger, see
\referto[Kuhna], p.\ 521) However, this result is true for {\it any}
(finite) cubical triangulation of $I^n$ which is related to the
fact that every open cover of $I^n$ is a uniform cover.
(Likewise every cubical triangulation of $I^n$
with rational vertices has
a subdivision which is regular in the above sense.)

On the
other hand, regular subdivisions are {\it necessary} when we deal
with infinite-dimensional cubes. Indeed, we will show
that the straightforward extension of Sperner's lemma
to arbitrary cubical subdivisions of $U_{\infty}$ is false.
Let $\c{V}$ be an open cover of $U_{\infty}$ such that
$\diam(V) < 1$ for each $V\in \c{V}$. Since $U_{\infty}$
is paracompact, we can assume that $\c{V}$ is locally
finite; let $\c{W}$ be an open refinement of $\c{V}$ such
that each $W\in\c{W}$  meets only finitely many members of $\c{V}$.
We can find a cubical subdivision $\c{K}$  of $U_{\infty}$ such that
$\c{K} \prec \c{W}$; i.\ e., for each cube $K\in\c{K}$ there
is $W\in\c{W}$ with $K\subset W$. Indeed, let $\c{K}_1$ be the
regular subdivision of $U_{\infty}$ into cubes of sidelength
$1/2$. Let $\c{K}'_1$ be the subset of all $K\in \c{K}_1$ such that
$K\subset W$ for some $W\in \c{W}$, and let $\c{K}'_2$ be the
cubical complex obtained by subdividing each
$K\in \c{K}_1\setminus \c{K}'_1$ into cubes of sidelength $1/4$.
Let $\c{K}_2$ denote the subset of all $K\in\c{K}_2$ such that
$K\subset W$ for some $W\in \c{W}$, and define $\c{K}_3$ as the
subdivision of $\c{K}_2\setminus \c{K}'_2$ into cubes of
sidelength $1/8$. Continue in this fashion ad infinitum. Then
let
$ \c{K}^* = \bigcup\{\c{K}'_n: n\in\omega\}.$
The set $\c{K}^*$ is naturally partially ordered with respect  to the
relation of inclusion. With this partial order, $\c{K}^*$ is a
tree. Furthermore, this tree is well-founded, i.\ e.\  each
maximal linearly ordered subset (that is, a decreasing sequence of cubes)
is finite. To see this, suppose that $K_1\supset K_2\supset\ldots$
is a decreasing sequence of cubes. Since
$\hbox{\rm diam}(K_n)\leq 2^{-n}$ and since $U_{\infty}$ is complete
as a metric space, there is $p\in \bigcap\{K_n: n\in \omega\}.$
Let $p\in W_p$, $W_p\in \c{W}$. Then $K_n\subset W_p$ already for
some $n$, contradicting the assumption that $K_{n+1}\in \c{K}^*$.
Thus, $\c{K}^*$ is well-founded. Consequently the minimal elements
of $\c{K}^*$ form the desired complex $\c{K}$. We define a
colouring $\varphi :\c{K}^{(0)}\to \{0,1\}^{\omega}$ as follows.
It is easy to see that the cardinality of the set $\c{K}$ is
$2^{\omega}$.
Choose for
each $p\in \c{K}^{(0)}$ some $V_p\in\c{V}$ such that $p\in V_p$,
and let $\c{V}' = \{V_p: p\in \c{K}^{(0)}\}$.
Then the cardinality of $\c{V}'$ is $2^{\omega}$ and
there is a bijection $\phi : \c{V}'\to \{0,1\}^{\omega}$.
Define $\varphi(p) = \phi(V_p)$. Then $\varphi$ satisfies
Sperner Condition since $\hbox{\rm diam}(V) < 1$ for each
$V\in\c{V}'$. However, for each cube $K\in \c{K}$ the
set $\varphi(K)$ of colours is finite, since $K$ meets
only finitely many members of $\c{V}$.

The above example also shows that the straightforward
extension of Lebesgue's Covering Dimension Theorem to the
infinite-dimensional setting is false: there is no $x\in U_{\infty}$
such that $\c{V}_x = \{V\in \c{V}: x\in V\}$ is infinite. Anyhow,
we prove that this extension is true for {\it uniform} covers
of $U_{\infty}$. We note here that if unit cubes of other
Banach spaces, e.\ g.\ $c_0(\omega)$ are considered, then the infinitary
Sperner lemma does not hold even for regular cubical subdivisions.
Indeed, $c_0(\omega)$ is separable and hence every uniform
cover of $c_0(\omega)$ has a uniformly locally finite
uniform refinement (\referto[Isbella], p.\ 111). Thus, by repeating the
construction in the above example, we can find arbitrarily
fine regular cubical subdivisions $\c{K}$ of the unit cube
of $c_0(\omega)$ with colourings $\varphi: \c{K}^{(0)}\to \omega$
satisfying Sperner Condition such that $\varphi(K^{(0)})$ is finite
for each $K\in \c{K}$. Combinatorially these subdivisions
correspond to the sets $[m]^{<\omega}$ of all sequences
$\sigma \in [m]^{\omega}$ satisfying $\sigma(i) = 0$ for
almost all $i\in\omega$.
It is, however, possible to give
modified versions of Sperner's lemma even in these cases;
we will return to this topic in Section 7.

\section{4. Definitions.} In this section we give definitions,
in addition to
those given in the previous section, necessary for the
proof of the main result. Let $n > 1$ be fixed.
The sum $\sigma + \tau$ of two elements $\sigma,\tau\in [n]^{\omega}$
is always understood relative to the interval $[n]$,
i.\ e., $(\sigma+\tau)(i) = \min(n,\sigma(i)+\tau(i))$
for all $i\in\Bbb N$.
We define for each $\sigma\in [n]^\omega$ the ``positive cube''
$K_\sigma$ with index $\sigma$ as the set of all
$\sigma+\tau$, where $\tau\in\{0,1\}^\omega$.
We also define a combinatorial
generalization of metric balls. Let $\sigma\in [n]^{\omega}$,
let $A\subset\omega$ and let $k\in [n]$. Then $B(\sigma,A,k)$
denotes the set of all $\tau\in [n]^{\omega}$ such that
$|\sigma(i) - \tau(i)| \leq k$ for $i\in A$ and $\sigma(i) = \tau(i)$
for $i\in\omega\setminus A$. In the proof of
5.1 we primarily consider those subsets $A$ of
$\omega$ for which both $A$ and $\omega\setminus A$ are infinite.
For any subset $S\subset \omega$ let $\c{A}(S)$ denote
the collection of all $A\subset \omega$ such that $S\subset A$ and
$|\omega\setminus A| = |A\setminus S| = \omega$. We define
$$ \hat B(\sigma, A, L) \; = \; \bigcup\{B(\sigma, A'\setminus A, L):
                A'\in \c{A}(A)\}.$$
Suppose that
$\c{G}$ is a covering of $[n]^\omega$.
We define here a
number that can be called a {\it local relative Lebesgue number} of the cover
$\c{G}$.
Given $\sigma \in [n]^\omega$ and
$A\subset\omega$, we define
$$\ell(\sigma,A,\c{G}) =
      \max\{k\in [n]: \exists G\in {\cal G}
                         (\hat B(\sigma,A,k)\subset G)\}.$$
We observe that $A_1\subset A_2$ implies
$\ell(\sigma,A_1,\c{G}) \leq \ell(\sigma,A_2,\c{G})$.
For each $\sigma\in
[n]^{\omega}$ there is $A\in\c{A}(\emptyset)$ such that
$\ell(\sigma,A,\varphi) = \ell(\sigma,A',\varphi)$ for all
$A'\in\c{A}(A)$. We also notice that if
$\sigma'\in \hat B(\sigma, A, k)$, say
$\sigma'\in B(\sigma, A'\setminus A, k)$, then
$\hat B(\sigma',A',k)\subset \hat B(\sigma, A, k)$, which
implies $\ell(\sigma', A', \c{G})\geq \ell(\sigma,A,\c{G)}$.

To facilitate the proof
of our first result, we define here a property $M$ that depends
on 5 parameters. (Unfortunately, simple arguments
such as that of \referto[Cohena] do not seem applicable
in this infinitary situation.)
Let $S$ be the set of all
$\sigma\in [n]^\omega$ such that $\sigma(i) = 0$ for
infinitely many $i\in \omega$.
Let $\sigma\in S$, let
$A\in\c{A}(\emptyset)$, let $k\in\omega$, let
$\c{G}$ be a covering of $[n]^\omega$,
and let $G\in \c{G}$. Then
$M(\sigma,A,k,G,\c{G})$ iff
$\hat B(\sigma,A,k)\subset G$ but
$\hat B(\sigma',A',k+1)\not\subset G'$ for all $G'\in \c{G}$ and 
for all extensions
$\sigma'$ of $\sigma$ in $\hat B(\sigma,A,k)$, where
$\sigma'\in S$ and $A'\in\c{A}(A)$.

Let $A\subset\omega$ and let $k\in\omega$. An $(A,k)$-function
is a function $\chi: \omega\to [-k,k]$ such that
$\chi(i) = 0$ for $i\in\omega\setminus A$. Thus, every
element of $B(\sigma,A,k)$ can be represented in the form
$\sigma + \chi$, where $\chi$ is an $(A,k)$-function.

\section{5. Infinitary Covering Dimension.} As
will be seen in the remarks following the proof
of 5.1 (Remark 5.4), the infinitary
Sperner lemma does not yield a direct extension of
Brouwer's Fixed Point Theorem. However, it is equivalent
to an infinitary extension
of Lebesgue's Covering Dimension Theorem. Here one has to
consider uniform covers instead of open covers, because the
cubes $U_{\infty}$ are not compact. A uniform space
$\mu X$ (see \referto[Isbella] for terminology) is called point-finite
if every uniform cover $\c{U}\in\mu$ has a uniform refinement
$\c{V}$ such that $\c{V}_x = \{ V\in \c{V}: x\in V\}$ is finite
for each $x\in X$.  Is every uniform (e.\ g., metric) space
point-finite? This question of A.\ H.\ Stone \referto[Stonea] and
J.\ Isbell \referto[Isbella] was answered in the negative by
E.\ V.\ \v S\v cepin \referto[Scepina] and
J.\ Pelant \referto[Pelanta]. \v S\v cepin
proved that $\ell_{\infty}(\beth_{\omega})$ is not point-finite,
where the beth number $\beth_{\omega}$ is defined inductively
by $\beth_0 = \omega$, $\beth_{n+1} = 2^{\beth_n}$ and
$\beth_{\omega} = \sup\{\beth_n: n\in\omega\}$. Pelant proved that
even $\ell_{\infty}(\beth_1)$ is not point-finite, by using his
combinatorial technique of cornets. By using graph-theoretic
results of Erd\"os, Galvin and Hajnal
\referto[Erdosa], V.\ R\"odl \referto[Rodla] has
given a simple proof showing that there is a non-point-finite
space of cardinality $\omega_1$. By using 5.1, we can easily
prove that $\ell_{\infty}(\omega)$ is not point-finite, by
showing that the subspace $U_{\infty}$ satisfies an
infinitary version of Lebesgue's Covering Dimension Theorem.
This result has also been announced by Pelant.

\th Theorem 5.1: Let $\c{U}$ be a uniform cover of
$U_{\infty}$ such that $\diam(U) < 1$ for each $U\in \c{U}$.
Then there is $x\in U_{\infty}$ such that $\c{U}_x$
is infinite.

\proof To facilitate the argument, we move from the covering
${\cal U}$ of $U_{\infty}$
to a covering ${\cal G}$ of $[n]^\omega$ for a suitable $n$.
Let $n\geq 2$ be such that the metric balls
$B(x,2/n)$ of $U_\infty$ refine ${\cal U}$. Then
for each $U\in {\cal U}$ let $G_U\subset[n]^\omega$ consist
of all $\sigma$ such that $(\sigma(i)/n)\in U$,
and define ${\cal G} = \{G_U : U\in {\cal U}\}$.
We shall construct a sequence of 4-tuples
$<\sigma_k,A_k,L_k,U_k>$, where $\sigma_k\in [n]^{\omega}$,
$A_k\in\c{A}(A_{k-1})$, $L_k\in [n]$, $U_k\in {\cal U}$,
such that $M(\sigma_k,A_k,L_k,G_{U_k},{\cal G})$ (as defined above)
holds for each $k$.

Let $A_0 = \emptyset$, let $S$ be as above and let
$$ L_1 = \max\{\ell(\sigma,A,{\cal G}):
              \sigma\in S, A\in\c{A}(A_0)\};$$
say $L_1 = \ell(\sigma_1,A_1,{\cal G})$, and let
$\sigma_1\in S$ and $U_1\in {\cal U}$
be such that $\hat B(\sigma_1,A_1,L_1)\subset G_{U_1}$.
Notice in particular that $L_1 < n$; this follows from the assumption
that $\diam(U) < 1$ for each $U\in \c{U}$.
We can assume that $\sigma_1(i) = 0$ for all $i\in \omega\setminus A_1$
and that there is an element $n_1\in A_1$ such that
$\sigma_1(n_1) = 0$.

It is clear that
$\hat B(\sigma, A, L_1+1)\not\subset G_{U}$ for all
$\sigma\in S$, $A\in \c{A}(A_1)$ and
$U\in {\cal U}$. It follows that
$M(\sigma_1,A_1,L_1,G_{U_1},{\cal G})$ holds.
For the inductive
hypothesis, assume that we have a sequence of 4-tuples
$<\sigma_k,A_k,L_k,U_k>$, $1\leq k\leq m$,
with the following properties:
\medskip
\item{1)} $M(\sigma_k,A_k,L_k,G_{U_k},{\cal G})$ holds for
          each $k\in \{1,\ldots,m\}$;
\smallskip
\item{2)} $L_1\geq \ldots \geq L_m$;
\smallskip
\item{3)} $A_{k+1}\in\c{A}(A_k)$ for each
          $k\in \{0,\ldots,m-1\}$;
\smallskip
\item{4)} there are fixed elements $n_i\in A_i$
   such that $\sigma_i(n_i) = 0$, where we assume that 
   $n_i \in A_i\setminus A_{i-1}$ for $i > 0$ and 
   that $|A_i\setminus A_{i-1}| > 1$;
\smallskip
\item{5)} if  $i\leq k\leq m$, then
          $\sigma_k\upharpoonright A_i =
                         \sigma_i\upharpoonright A_i;$
\smallskip
\item{6)} if $i\in \omega\setminus A_m$,
          then $\sigma_k(i) = 0$ for all $k\in\{1,\ldots,m\}$;
\smallskip
\item{7)} $\hat B(\sigma_{i+1}, A_{i+1}, 1)\not\subset G_{U_i}$
          for all $i\in \{1,\ldots,m-1\}$.
\medskip
We shall construct a 4-tuple
$<\sigma_{m+1},A_{k+1},L_{k+1},U_{k+1}>$ such that the above
conditions 1) - 7) hold with $m$ replaced by $m+1$.
Notice that for each
$A\in \c{A}(A_m)$ and each $(A\setminus A_m,L_m)$-function
$\chi:\omega\to [-L_m,L_m]$, one has
$\sigma_m + \chi \in G_{U_m}$.
We claim that there is
$A'_m\in\c{A}(A_m)$ and an $(A'_m\setminus A_m,L_m)$-function
$\chi_m$ such that
$$ B(\sigma_m+\chi_m,A'_m\setminus A_m,1)\not\subset G_{U_m}.$$
Indeed, suppose that there is no such $A'_m$. Then
$\hat B(\sigma_m, A_m, L_m +1)\subset G_{U_m}$. To see this,
let $\alpha\in \hat B(\sigma_m,A_m,L_m+1)$. Thus, there is
$A\in \c{A}(A_m)$ such that
$\alpha\in B(\sigma_m,A',L_m+1)$, where $A' = A\setminus A_m$.
We have
$|\sigma_m(i) - \alpha(i)|\leq L_m+1$ for all $i\in A'$
and $\alpha(i) = \sigma_m(i)$ for $i\in\omega\setminus A'$.
Define a function $\beta\in [n]^{\omega}$ by setting
$\beta(i) = \sigma_m(i)$ for $i\in \omega\setminus A'$, set
$\beta(i) = \alpha(i)$ for $i\in A'$ such that
$|\alpha(i) - \sigma_m(i)|\leq L_m$ and otherwise
$\beta(i) = \sigma_m(i) - L_m$ or $\beta(i) = \sigma_m(i) + L_m$
depending on whether $\alpha(i) < \sigma_m(i)$ or
$\alpha(i) > \sigma_m(i)$. Clearly
$\beta\in B(\sigma_m, A', L_m)$ and $||\alpha - \beta||\leq 1$,
and thus $\alpha\in B(\sigma_m + \chi, A', 1)$ for the
$(A',L_m)$-function  $\chi = \beta - \sigma_m$. Therefore,
by our assumption,
we have $\alpha\in G_{U_m}$ and consequently
this shows that $\hat B(\sigma_m, A_m, L_m+1)\subset G_{U_m}$,
which is a contradiction with the definition of $L_m$.
({\it Notice that for this contradiction we need the crucial
property that $L_m < n$.})
Thus, the desired function $\chi_m$ and the desired
set $A'_m\in\c{A}(A_m)$ exist.

Let
$$ L_{m+1} = \max\{\ell(\sigma, A, {\cal G}): \sigma\in E, A\in\c{A}(A'_m)\},$$
where $E$ denotes the set of all extensions of $\sigma_m+\chi_m$
in $\hat B(\sigma_m, A_m, L_m)$.
It is easy to see that $L_{m+1}\leq L_m$. We can find
$\sigma_{m+1}$, $A_{m+1}$, $U_{m+1}$ with the following
properties:
\medskip
\item{1)} $A_{m+1}\in\c{A}(A'_m)$;
\smallskip
\item{2)} $\hat B(\sigma_{m+1}, A_{m+1}, L_{m+1})\subset G_{U_{m+1}};$
\smallskip
\item{3)} $\hat B(\alpha, A, L_{m+1}+1)\not\subset G_{U}$
          for all $A\in \c{A}(A_{m+1})$, all $U\in \c{U}$
          and all extensions $\alpha\in S$ of $\sigma_{m+1}$ in
          $\hat B(\sigma_{m+1}, A_{m+1}, L_{m+1})$.
\medskip
\noindent
%
%
It follows that
$M(\sigma_{m+1}, A_{m+1}, L_{m+1}, G_{U_{m+1}},\c{G})$ holds.
Finally, we note that Condition 4) can easily be satisfied
since $A_{m+1}$ can be replaced by a larger infinite set.
This finishes the inductive step.

As $L_{i+1}\leq L_i$ for $i\in \Bbb N^*$, there is
(the least) $i_0\in \Bbb N$
such that $i\geq i_0$ implies $L_i = L_{i_0}$.
(Moreover, notice that $L_i\geq 1$ for all $i$ by the
choice of $n$.)
Define
$$ \hat\sigma = \lim \sigma_i,$$
i.\ e.,    $\hat\sigma(i) = \sigma_k(i)$ for $i\in A_k$ and
$\hat\sigma(i) = 0$ otherwise. Then
$\hat\sigma\in \hat B(\sigma_i,A_i, L_{i_0})$ for all $i\geq i_0$.
Indeed, the support of $\hat\sigma$ is contained in $A_{\omega} =
\bigcup\{A_k: {k\in\Bbb N}\}\setminus\{n_k: {k\in \Bbb N^*}\}$,
and this is -- by the inductive
construction --  an element of $\c{A}(\emptyset)$. For each
$i\geq i_0$, we have $\hat \sigma = \sigma_i + \chi_i$, where
$\chi_i$ is an $(A,L_{i_0})$-function with $A\in \c{A}(A_i)$.
We claim that $i,j\geq i_0$, $i\neq j$ implies $U_i\neq U_j$.
To prove this, let us assume $i<j$ and $U_i = U_j$ to derive
a contradiction.
By our assumption and
by the choice of the sets $G_{U}$, we have
$$\hat B(\sigma_j, A_j, L_i)\subset G_{U_j}=G_{U_i},\leqno(*)$$
and therefore
$B(\sigma_i + \chi_i, A'_i\setminus A_i, 1)\subset G_{U_i}$.
(Recall the above definition of the set $A'_i$.) To prove this claim,
suppose that $\xi\in B(\sigma_i + \chi_i, A'_i\setminus A_i, 1)$.
For
$k\in A'_i\setminus A_i$ we have
by the definition of $\sigma_{i+1}$ that
$\sigma_j(k) = \sigma_i(k) + \chi_i(k)$.
Notice that $\xi(k) = 0$ for all
$k\in \omega\setminus A'_i$.
If $k\in A_j\setminus A'_i$, then
$|\sigma_i(k) - \sigma_j(k)| \leq L_{i_0}$. Finally, if
$k\in A_i$, then $\sigma_j(k) = \sigma_i(k)$. It then follows
from $L_{i_0}\geq 1$ that
$|\xi(k) - \sigma_j(k)| \leq L_{i_0}$ for all
$k\in A_j$, and hence by $(*)$
we have $\xi\in G_{U_i}$. This contradiction proves that
$U_i\neq U_j$.

Finally, we have
$|(\c{G})_{\hat\sigma}|\geq\omega$.
In fact, as
$||\sigma_i - \hat\sigma||\leq L_{i_0}$ for all
$i\geq i_0$ and by the inductive construction
$\hat B(\sigma_i, A_i, L_{i_0})\subset G_{U_i}$, we have
$\hat\sigma\in G_{U_i}$ for all $i\geq i_0$. But then
$x\in U_i$ for infinitely many $i$, where
$x = (\hat\sigma(i)/n)$. This concludes the proof of 5.1.

\noindent
{\bf Remark 5.2}: The statement of 5.1
is the best possible. One cannot prove consistently with
ZFC that under the hypotheses of 5.1 there is $x\in {U_\infty}$
such that $|({\cal U})_x| = \kappa > \omega$. Indeed,
assume that {\bf CH} (the continuum hypothesis) holds. Then
the uniformity of
$U_{\infty}$ has a basis of covers of cardinality $\omega_1$.
But then the proof given for the point-finiteness of separable
spaces in \referto[Isbella], p.\ 111, shows that the uniformity of
$U_{\infty}$ has a basis
consisting of point-countable covers. Moreover, by using the method
of Section 3 we can thus construct a regular cubical triangulation
$\c{K}$ of $U_{\infty}$ and a Sperner colouring
$\varphi: \c{K}^{(0)}\to \{0,1\}^{\omega}$ such that
$\varphi(K)$ is at most a countable set for
each cube $K$ of $\c{K}$.

\noindent
{\bf Remark 5.3}: The simplest proof showing that a metric
(uniform) space is not point-finite is given by Pelant and R\"odl
in \referto[Pelantc]. In fact, they implicitly
formulate and prove a ``weak'' infinitary
Sperner theorem. Suppose that $m,n\in\omega$, $n > 0$, and let
$\varphi: [\beth_{m+n-1}]^n\to\beth_m$ be a mapping
(``colouring'') such that $a\cap a'=\emptyset$
implies $\varphi(a)\neq\varphi(a')$ for all
$a,a'\in[\beth_{m+n+1}]^n$ (``Sperner Condition''). Then there
exists a subset (``simplex'') $\Delta\subset [\beth_{m+n-1}]^n$
such that 1) $|\Delta| = \beth_m$; 2) $a\neq a'$ implies
$\varphi(a) \neq\varphi(a')$ for all $a,a'\in \Delta$ and
3) $|\bigcap\Delta| = n-1$. The proof (by induction on $n$)
easily follows from the strong assumptions. This result is used
to show that $\ell_1(\beth_{\omega})$  is not point-finite.

\noindent
{\bf Remark 5.4}: Theorem 5.1 as such does not imply a useful fixed
point theorem for mappings $U_{\infty}\to U_{\infty}$. Indeed, the
regular cubical triangulations correspond to uniformly continuous
mappings $f :U_{\infty}\to U_{\infty}$, but the usual method
of using Sperner's lemma (see e.\ g.\ \referto[Knastera]) only yields
that for each $\epsilon > 0$ there is an infinite set
$A\subset\omega$ and $x\in U_{\infty}$ such that
$|x_i - (f(x))_i| < \epsilon$ for $i\in A$. This can readily
be proved without the use of 5.1. It has been shown
(\referto[Benyaminia])
that there are even Lipschitz mappings $f: U_{\infty}\to U_{\infty}$
without approximate fixed point; i.\ e.,  there is $\epsilon > 0$
such that $||x - f(x)||_{\ell_{\infty}}\geq \epsilon$ for all
$x\in U_{\infty}$. This leads us to the following question.

\noindent
{\bf Question:} Let $f: U_{\infty}\to U_{\infty}$ be a uniformly
continuous mapping. Is there an infinite subset $A\subset\omega$
such that $(f(x))_i = x_i$ for $i\in A$?

\section{6. An Infinitary version of Sperner's Lemma.} In this section
we state and prove
our infinitary version of Sperner's lemma. Let us recall
that a mapping $\varphi :[n]^{\omega}\to \{0,1\}^{\omega}$ is a
Sperner colouring if $\varphi(\sigma)\neq\varphi(\sigma')$
whenever $||\sigma - \sigma'|| = n$.

\th Theorem 6.1: Let $n\in \Bbb N$ and let
$\varphi :[n]^{\omega}\to \{0,1\}^{\omega}$ be a Sperner colouring.
Then there is $\sigma\in [n]^{\omega}$ such that
$\varphi(K_{\sigma})$ is infinite.

\proof
We will define a uniform cover $\c{U}$ of
$U_{\infty}$ and apply
Theorem 5.1. For each $\sigma\in [n]^\omega$ define the
set $G_\sigma$ as the product
$$ G_\sigma = \prod_{k\in\Bbb N} I_k(\sigma),$$
where for each $k\in\Bbb N$,
$I_k(\sigma)$ is the open interval
$]{\sigma(k)-2/3\over n},{\sigma(k)+2/3\over n}[$ for  $0 < \sigma(k) < n$,
the interval $[0,1/n[$ for $\sigma(k) = 0$
and $]1-1/n, 1]$ for $\sigma(k) = n$.
For each $\tau\in \{0,1\}^\omega$, let
$$ U_\tau = \bigcup\{G_\sigma : \varphi(\sigma) = \tau\}.$$
Then $\c{U} = \{U_\tau : \tau\in \{0,1\}^\omega\}$ is
a uniform (open) cover of $U_\infty$, and
$\diam(U_\tau) < 1$ for all $\tau$, because
$\varphi$ is a Sperner colouring.
By 5.1 there is $x\in U_\infty$ such that
$(\c{U})_x$ is infinite. Thus,
$x$ is contained in infinitely many
sets $G_\sigma$, each mapped by
$\varphi$ to a different $\tau$. Let
$\Sigma$ be an infinite set of elements
$\sigma$ such that $x\in G_\sigma$ and
which (pairwise) map to distinct colours.
Given two such elements $\sigma_1,\sigma_2$,
we have $|\sigma_1(i) - \sigma_2(i)| \leq 1$ for
all $i$, by the definition of the sets $G_\sigma$.
It follows that there is a cube $K_\sigma$ for which
$\Sigma$ forms a subset of vertices; indeed,
we may define $\sigma$ as the coordinatewise
infimum of the elements of $\Sigma$. This
proves 6.1.

\noindent
{\bf Remark 6.2}: In the same way as the classical
Sperner lemma corresponds to a homology theory of simplicial
complexes (see, e.\ g.\  \referto[Browna]),
our infinitary version of Sperner's lemma
(Theorem 6.1) corresponds to an infinitary homology theory
of infinite-dimensional cubical complexes.

\section{7. The case of $c_0$.}
As noted earlier, the regular cubical triangulations of the
unit ball of the Banach space $c_0$ correspond to the sets
$[n]^{<\omega}$ of all $\sigma\in [n]^{\omega}$ such that
$s(\sigma) = \{k\in\omega: \sigma(k)\neq 0\}$ is finite.
We also noted that the infinitary version of Sperner's
lemma does not hold for these sets. However, although
there are Sperner colourings $\varphi:[n]^{<\omega}\to \omega$
such that for each cube $K_{\sigma}$ the vertices
are coloured with only finitely many colours, one can
show that there is a sequence $(K_{\sigma_k})$
of cubes and a sequence
$(\tau_k)_{k\in\Bbb N}$ of distinct colours $\tau_k\in\omega$
such that
\medskip
\item{1)} $\{\tau_1,\ldots,\tau_k\}\subset \varphi(K_{\sigma_k})$;
\smallskip
\item{2)} $\sigma_{k+1}$ is an extension of $\sigma_k$ for each
$k$, i.\ e., $\sigma_{k+1}(i) = \sigma_k(i)$ for
$i\in s(\sigma_k)$, where $s(\sigma)$ denotes the support of $\sigma$.
\medskip
\noindent
This result is obtained from the following version of
Lebesgue's Covering Dimension Theorem for the unit
cube $U(c_0)$ of $c_0$.
It is established by virtually the same proof as that
given for 5.1 except that the families $\c{A}(S)$ are
replaced by the families $\c{F}(S)$ of finite subsets.
(Let us note that even this result was announced
by Pelant in 1986. The proof was based on his technique
of cornets.)

\th Theorem 7.1: Let $\c{U}$ be a uniform covering
of $U(c_0)$ such that $\diam(U)< 1$ for each $U\in\c{U}$.
Then there is a sequence $(U_n)_{n\in\Bbb N}$ of elements
of $\c{U}$ such that for each $n\in\Bbb N$, we have
$U_1\cap\cdots\cap U_n\neq\emptyset$.

We will interpret 7.1 with respect to {\it Noetherian}
covers of uniform spaces. Let $X$ be a set, and let
$\c{V}$ be any point-finite cover of $X$. There is a
natural partially ordered set $\c{P}(\c{V})$ associated
with $\c{V}$ which consists of all finite subsets
$\{V_1,\ldots,V_n\}$ such that $V_1\cap\cdots\cap V_n\neq\emptyset$
and which is ordered with respect to set inclusion.
(The poset $\c{P}(\c{V})$ corresponds to a simplicial
complex in which the finite intersecting subsets are
regarded as simplices.) The cover $\c{V}$ is called
Noetherian if $\c{P}(\c{V})$ does not contain any
infinite increasing chain. Theorem 7.1 implies that
no bounded uniform cover of $U(c_0)$ is Noetherian.
Since the Lebesgue covering dimension of an $n$-cube
can be regarded as the minimum of the maximal length
of chains in posets $\c{P}(\c{V})$, where $\c{V}$
is a bounded open covering of the cube, Theorem 7.1
can again be considered an infinitary version of
Lebesque's Covering Dimension Theorem. As in the case of
$U(\ell_\infty)$, the extension fails for general open
coverings. Indeed, any paracompact space has a base
of open Noetherian coverings (this result has been
established independently by J.\ Fried \referto[Frieda]).

The infinite chains of intersecting sets are
representatives of infinite-dimensional simplices of
dimension $\omega$, and 7.1 corresponds to an
infinitary homology theory in the same way as the classical
Sperner lemma is related to the classical simplicial homology
theory. In this context, the cube $U(c_0)$
represents the {\it finitary boundary} of $U(\ell_\infty)$, to be compared
with $S^{n-1}$ as the boundary of $I^n$. These problems
will be considered in another paper.

{\bf Acknowledgement.} The author expresses his gratitude, once again, to Heikki 
Junnila for his helpful comments.

\bigskip
\centerline{\smallcaps References}
\bigskip
{
\beginref

\ref Benyamini, Y., and Y.\ Sternfeld: Spheres in infinite-dimensional
     normed spaces are Lipschitz contractible.- Proc.\ Amer.\ Math.\
     Soc. 88:3, 1983, pp.\ 439--445.[Benyaminia]

\ref Brown, A.\ B., and S.\ Cairns: Strengthening of Sperner's lemma
     applied to homology theory.- Proc.\ Nat.\ Acad.\ Sci.\ U.\ S.\ A.\
     47, 1961, pp.\ 113--114.[Browna]

\ref Cohen, D.\ I.\ A: On the Sperner Lemma.- Journal of Comb.\
     Theory 2, 1967, pp.\ 585--587.[Cohena]

\ref Engelking, R.: Dimension theory.- Polish Scientific Publishers,
     Warszawa, 1978.[Engelkinga]

\ref Erd\"os, P., Galvin, F., and A.\ Hajnal: On set-systems having
     large chromatic number and not containing prescribed
     subsystems.- Infinite and Finite Sets (A.\ Hajnal, R.\ Rado,
     V.\ T.\ S\'os (Eds.)), North-Holland, 1976, pp.\ 425--513.
     [Erdosa]

\ref Fan, Ky.: A generalization of Tucker's combinatorial lemma
     with topological applications.- Ann.\ of Math. 2, 56, 1952,
     pp.\   431--437.[Fana]

\ref Fried, J: Personal communication.-{}[Frieda]

\ref Goebel, K.: On the minimal displacement of points under
     Lipschitzian mappings.- Pacific J.\ Math. 45, 1973, pp.\ 151--163.
     [Goebela]

\ref Isbell, J.\ R.: Uniform spaces.- Mathematical Surveys 12,
     Amer.\ Math.\ Soc., Providence, Rhode Island, 1964.[Isbella]

\ref Knaster, B., Kuratowski, C.\ and S.\ Mazurkiewicz.: Ein
     Beweis  des Fixpunksatzes f\"ur $n$-dimensionale Simplexe.-
     Fund.\ Math. 14, 1929, pp.\ 132--137.[Knastera]

\ref Kry\'nski, S.:  Remarks on matroids and Sperner's lemma.-
     Europ.\ J.\ Combinatorics 11, 1990, pp.\ 485--488.[Krynskia]

\ref Kuhn, H.\ W.: Some combinatorial lemmas in topology.-
     IBM J.\ Res.\ and Dev. 4, 1960, pp.\ 518--524.[Kuhna]

\ref Lindstr\"om, S.: On matroids and Sperner's lemma.-
     Europ.\ J.\ Combinatorics 2, 1981, pp.\ 65--66.[Lindstroma]

\ref L\'ovasz, L.: Matroids and Sperner's lemma.- Europ.\ J.\
     Combinatorics 1, 1980, pp.\ 65--66.[Lovasza]

\ref Mani, P.: Zwei kombinatorisch-geometrische S\"atze vom
     Typus Sperner-Tucker-Ky Fan.- Mo\-nats\-hefte f.\ Ma\-the\-ma\-tik 71,
     1967, pp.\ 427--435.[Mania]

\ref Pelant, J.: Combinatorial properties of uniformities.- General
     Topology and its Relations to Modern Analysis and Algebra
     IV, Lecture Notes in Mathematics 609, Springer-Verlag,
     Berlin-Heidelberg-New York, 1977, pp. \ 154--165.[Pelanta]

\ref Pelant, J.: Embeddings into $c^+_0(\omega)$.- preprint.
     [Pelantb]

\ref Pelant, J., and V.\ R\"odl: On coverings of infinite-dimensional
     metric spaces.- preprint.[Pelantc]

\ref R\"odl, V.: Small spaces with large point-character.-
     Europ.\ J.\ Combinatorics 8, 1987, pp.\ 55--58.[Rodla]

\ref Sperner, E.: Neuer Beweis f\"ur die Invarianz der Dimensionzahl
     und des Gebietes.- Abh.\ Math.\ Sem.\ Hamburg 6, 1928,
     pp.\ 265--272.[Spernera]

\ref Sperner, E.: Kombinatorik bewerter Komplexe.- Abh.\ Math.\ Sem.\
     Hamburg 39, 1973, pp.\ 21--43.[Spernerb]

\ref Stone, A.\ H.: Universal spaces for some metrizable uniformities.-
     Quart.\ J.\ Math.\ Oxford, Ser. 11, 1960, pp.\ 105--115.[Stonea]

\ref \v S\v cepin, E.\ V.: On a problem of Isbell.- Soviet Math.\
     Dokl.\ 16, 1975, pp.\ 685--687.[Scepina]

\ref Tucker, A.\ W.: Some topological properties of disk and
     sphere.- Proc.\ First Canadian Math.\ Congress,
     Montreal, Canada, 1945, pp.\ 285--309.[Tuckera]

\endref
}

\vskip2cm
\section{The address:}
\bigskip
\hbox{\vbox{
            \hbox{Aarno Hohti}
            \hbox{University of Helsinki}
            \hbox{Department of Mathematics}
            \hbox{Yliopistonkatu 5}
            \hbox{SF-00100 Helsinki}
            \hbox{FINLAND}} \hfill}

\end